\input amssym.def
\input amssym.tex
\def\C{{\Bbb C}}
\def\R{{\Bbb R}}
\def\Q{{\Bbb Q}}
\def\qed{{~~~\vrule height .75em width .4em}}
\def\cent#1#2{{\bf C}_{#1}(#2)}
\def\deg#1{{\rm deg}(#1)}
\def\nor{\triangleleft\,}
\def\gen#1{\langle#1\rangle}
\def\snor{\mathrel{\hbox{$\triangleleft\triangleleft$}}}
\def\Oh#1#2{{\bf O}^{#1}(#2)}
\def\norm#1#2{{\bf N}_{#1}(#2)}
\def\ritem#1{\item{{\rm #1}}}
\def\iitem#1{\goodbreak\par\noindent{\bf #1}}
\def\inv{^{-1}}
\def\aut#1{{\rm Aut}(#1)}
\def\gal#1#2{{\rm Gal}(#1/#2)}
\def\sbs{\subseteq}
\def\dimm#1#2{{\rm dim}_{#1}(#2)}
\def\sps{\supseteq}
\def\minp#1#2{{\rm min}_{#1}(#2)}

\def\ref#1{{\bf[#1]}}
\def\book{2}
\def\monthly{1}
\def\XX{{\cal X}}
\def\YY{{\cal Y}}
\def\PP{{\cal P}}

\magnification = \magstep1
{\nopagenumbers
\font\bb = cmbx12
\vglue 1.5 truein
\centerline{{\bb REAL FIELDS AND}}
\smallskip
\centerline{{\bb REPEATED RADICAL EXTENSIONS}}
\medskip
\centerline{by}
\bigskip
\centerline{{\bf I.\ M.\ Isaacs}}
\centerline{and}
\centerline{{\bf David Petrie Moulton}}
\bigskip\bigskip
\centerline{Mathematics Department}
\centerline{University of Wisconsin}
\centerline{480 Lincoln Drive}
\centerline{Madison WI~~~53706}
\centerline{USA}
\bigskip
\centerline{E-mail:}
\centerline{isaacs@math.wisc.edu}
\centerline{moulton@math.wisc.edu}
\vfil
\eject}

\iitem{1. Introduction.}

Recall that a field extension $F \sbs L$ is said to be a {\bf radical
extension} if it is possible to write $L = F[\alpha]$, where $\alpha \in L$ is
an element with $\alpha^n \in F$ for some positive integer $n$. More
generally, an extension $F \sbs L$ is a {\bf repeated radical extension} if
there exist intermediate fields $L_i$ with $F = L_0 \sbs L_1 \sbs \cdots
\sbs L_r = L$ and such that each field $L_i$ is a radical extension of
$L_{i-1}$ for $0 < i \le r$.

Given a polynomial $f(X)$ over a field $F$ of characteristic zero, let $S$
be a splitting field over $F$ for $f$. Then as usual, we say that $f$ is
{\bf solvable by radicals} if $S$ is contained in some repeated radical
extension of $F$. A celebrated theorem of Galois asserts that this occurs
if and only if the associated Galois group $\gal SF$ is a solvable group.

It is well known that intermediate fields of repeated radical extensions
need not themselves be repeated radical extensions of the ground field.
The solvability of $\gal SF$, therefore, does not guarantee that $S$ is a
repeated radical extension of $F$, and so the phrase ``contained in" in
the statement of Galois' theorem is essential. For example, take
$F = \Q$, the rational numbers, and consider the polynomial
$f(X) = X^3 - 6X + 2$. It is easy to see that $f$ has three real roots, and
so we can take $S \sbs \R$. Of course, the cubic polynomial $f$ is
solvable by radicals; we can see this explicitly by calculating that the
three roots of $f$ are given by the formula $r = \alpha + 2/\alpha$, where
$\alpha$ runs over the three complex cube roots of the complex number
$-1+\sqrt{-7}$. If $S$ were a repeated radical extension of $\Q$, there
would have to be some alternative way to express these roots in terms of
{\it real} radicals. This is impossible, however, since $f$ is easily seen
to be irreducible, and it is a classical result that if an irreducible cubic
polynomial has three real roots, then these roots definitely are 
not expressible in terms of real radicals. More generally, we have the
following (known) result. (See \ref\monthly\ or Theorem~22.11 of
\ref\book. Also, we include a somewhat simplified proof here, in
Section~4.)
\medskip

\iitem{THEOREM A.}~~{\sl Let $Q$ be any subfield of the real numbers
$\R$ and suppose that $f \in Q[X]$ is irreducible and splits over $\R$. If
any one of the roots of $f$ lies in a real repeated radical extension of
$Q$, then $\deg f$ must be a power of $2$.}
\medskip

One of the principal results of this paper is the following, which shows
that at least in certain cases, intermediate fields of repeated radical
extensions are themselves repeated radical extensions. As in
Theorem~A, the ground field need not be the rational numbers; any field
$Q \sbs \R$ will suffice.
\medskip

\iitem{THEOREM B.}~~{\sl Suppose that $Q$ is a real field and that
$Q \sbs L$ is a repeated radical extension with $|L:Q|$ odd. If
$Q \sbs K \sbs L$, then $K$ is a repeated radical extension of $Q$.}
\medskip

We will show by example that the condition in Theorem~B that the
ground field $Q$ should be real cannot be dropped. Our proof of
Theorem~B begins by observing that it is no loss to assume that
$L \sbs \R$. To handle that case, we derive a useful, but somewhat
technical, characterization of real repeated radical extensions. This
characterization also has other applications and, in particular, it can be
used to prove the following result, which in some sense complements
Theorem~A.
\medskip

\iitem{THEOREM C.}~~{\sl Suppose that $Q$ is a real field and that
$f \in Q[X]$ is irreducible of odd degree. If $f$ has some root $\alpha$
in a real repeated radical extension of $Q$, then $\alpha$ is the only real
root of $f$.}
\medskip

In Theorem~B, we considered intermediate fields of odd-degree
repeated radical extensions. It is perhaps somewhat surprising that we
get a similar result in exactly the opposite case, where the degree of the
extension is a power of $2$. As we shall see, the derivation of this result
from our characterization of real repeated radical extensions is
considerably easier than is the proof of Theorem~B. 
\medskip

\iitem{THEOREM D.}~~{\sl Suppose that $Q \sbs L$, where $L$ is a real
repeated radical extension of $Q$. If $|L:Q|$ is a power of $2$ and
$Q \sbs K \sbs L$, then $K$ is a repeated radical extension of $Q$.}
\medskip

The full strength of the hypothesis that the fields are real is not actually
needed for Theorems~A and~D. For these results, it suffices that the
relevant field $L$ is {\bf quasireal}, which we define to mean that $L$ has
characteristic zero and that the only roots of unity it contains are $\pm1$.
In fact, the characterization of real repeated radical extensions to
which we referred earlier works more generally for quasireal fields.
In proving Theorems~B and~C, however, we shall use the realness
assumptions more fully, although even for those results, it would be
sufficient to assume that the fields are merely formally real; they need
not actually be subfields of $\R$. (Recall that a field is {\bf formally real}
if $-1$ is not a sum of squares, or equivalently, if the field can be
ordered.)

We close this introduction by acknowledging the contribution of
H.~W.~Lenstra,~Jr., who first suggested Theorem~B as a problem. Part
of this paper is part of the Ph.D.\ thesis of the second author, written
under Lenstra's supervision at the University of California, Berkeley.
Also, the first author is grateful to the University of California for its
hospitality while he was visiting there, on sabbatical from the University
of Wisconsin, Madison.
\bigskip

\iitem{2. Prime-degree extensions.}

We begin with an easy, but useful, lemma and some simple
consequences. Most of this material is well known.
\medskip

\iitem{(2.1) LEMMA.}~~{\sl Let $Q \sbs L$ be fields and suppose that
$L = Q[\alpha]$, where $\alpha^n$ lies in $Q$ for some integer $n \ge 1$.
Write $d = |L:Q|$. The following then hold.
\smallskip
\ritem{(a)} In any extension field of $L$, all roots of the minimal
polynomial $\minp Q\alpha$ have the form $\delta\alpha$, where
$\delta$ is an $n\,$th root of unity.
\smallskip
\ritem{(b)} We have $d \le n$, and if $\alpha^d \in Q$, then $d$
divides $n$.
\smallskip
\ritem{(c)} For some $n\,$th root of unity $\epsilon \in L$, we have
$\epsilon\alpha^d \in Q$. In particular, if $Q$ contains all $n\,$th roots
of unity in $L$, then $\alpha^d \in Q$.
\medskip}

\iitem{Proof.}~~Write $a = \alpha^n$, so that $\alpha$ is a root of
$X^n - a \in Q[X]$. The minimal polynomial $f = \minp Q\alpha$ must
therefore divide $X^n - a$, and hence each root $\beta$ of $f$ is also
a root of $X^n - a$. Thus $\beta^n = a$, and it follows that
$\beta = \delta\alpha$, as claimed.

Since $d = \deg f$, it follows that $\alpha^r$ cannot lie in $Q$ for any
positive exponent $r < d$, and in particular we have $d \le n$.
Writing $n = qd + r$ with $0 \le r < d$, we see that
$\alpha^r = \alpha^n(\alpha^d)^{-q}$, and this lies in $Q$ if
$\alpha^d \in Q$. It follows that $r$ cannot be positive in this
case, and thus $d$ divides $n$.

Also, as $\deg f = d$, it follows from (a) that the product of the $d$
roots of $f$ in a splitting field (counting multiplicities) has the form
$\epsilon\alpha^d$ for some $n\,$th root of unity $\epsilon$ in the
splitting field. But this product equals $\pm f(0)$, and so it lies in
$Q$, and hence in $L$. Since $\alpha^d \in L$, we deduce that
$\epsilon \in L$, as required.\qed
\medskip

Although it is not needed for what follows, we cannot resist mentioning
the following pleasant application of Lemma~2.1.
\medskip

\iitem{(2.2) COROLLARY}~~{\sl Let $f(X) = X^p - a \in Q[X]$, where
$Q$ is any field and $p$ is a prime number. Then either $f$ is irreducible,
or else it has a root in $Q$.}
\medskip

\iitem{Proof.}~~Let $\alpha$ be a root of $f$ in some extension field
$E = Q[\alpha]$ and write $m = |E:Q|$. If $m = p$, then $f$ is irreducible,
and so we assume that $m < p$. In particular, $m$ and $p$ are
coprime, and we can thus choose integers $k$ and $l$ such that
$mk + pl = 1$.

Since $\alpha^p \in Q$, we know by Lemma~2.1(c) that
$\epsilon\alpha^m \in Q$ for some $p\,$th root of unity $\epsilon$. Thus
$$
\epsilon^k\alpha = \epsilon^k\alpha^{mk}\alpha^{pl} =
(\epsilon\alpha^m)^k(\alpha^p)^l \in Q \, .
$$
Since $\epsilon^k\alpha$ is a root of $f$, this completes the proof.\qed
\medskip

\iitem{(2.3) LEMMA.}~~{\sl Let $Q \sbs L$ be a radical extension
and assume that $|L:Q| = p$, where $p$ is an odd prime.
\smallskip
\ritem{(a)} If $L$ is Galois over $Q$, then $L$ contains some root of unity
different from $\pm1$, and so $L$ is not quasireal.
\smallskip
\ritem{(b)} If $L$ is not Galois over $Q$, then $L = Q[\alpha]$ for some
element $\alpha$ with $\alpha^p \in Q$.
\medskip}

\iitem{Proof.}~~Write $L = Q[\alpha]$, where some power of $\alpha$
lies in $Q$, and consider the minimal polynomial $f = \minp Q\alpha$. If
$L$ is Galois over $Q$, then $f$ has at least $\deg f = p \ge 3$ distinct
roots in $L$. By Lemma~2.1(a), each of these roots has the form
$\delta\alpha$ for some root of unity $\delta \in L$, and thus $L$ contains
three different roots of unity, at least one of which must be different from
$\pm1$. This proves (a).

By Lemma~2.1(c), we know that $\epsilon\alpha^p \in Q$ for some root
of unity $\epsilon \in L$, and we have $Q \sbs Q[\epsilon] \sbs L$. If $L$
is not Galois over $Q$, however, then $L \ne Q[\epsilon]$. Since $|L:Q|$
is prime, we conclude that $Q[\epsilon] = Q$ and $\epsilon \in Q$, and it
follows that $\alpha^p \in Q$, as required.\qed 
\medskip

As is suggested by these lemmas, we shall need to control the roots
of unity in a field extension $Q \sbs E$. For this purpose, it will be
useful to consider intermediate fields $F$ that are abelian extensions of
$Q$ and contain all roots of unity in $E$. (Recall that a field extension
$Q \sbs F$ is said to be {\bf abelian} if it is a Galois extension such that
$\gal FQ$ is an abelian group. Also, we include in the definition
of ``Galois" the assumption that the extension has finite degree.)
Given any finite degree extension $Q \sbs E$, it is clearly always possible
to find such a field $F$: simply take $F$ to  be the field generated by $Q$
and all roots of unity in $E$. It is useful to have a little more freedom in
selecting $F$, however, and this is the purpose of the following result.
\medskip

\iitem{(2.4) LEMMA.}~~{\sl Let $Q \sbs K \sbs E$, where $K$ is
abelian over $Q$. Then there exists a field $F$ with $K \sbs F \sbs E$
such that $F$ is abelian over $Q$ and contains all roots of unity in $E$.}
\medskip

\iitem{Proof.}~~Let $L$ be the field generated over $Q$ by all roots of
unity in $E$. Then $L$ is Galois over $Q$, and hence the compositum
$F = KL$ of $K$ and $L$ in $E$ is also Galois over $Q$.
Furthermore, since each of $L$ and $K$ is abelian over $Q$, it is easy to
see that $F$ is also abelian over $Q$, and this completes the proof.\qed
\medskip

Given a characteristic zero field extension $Q \sbs L$, we want to be
able to determine whether or not it is a repeated radical extension. We
begin by choosing an arbitrary Galois extension $E$ of $Q$ that
contains $L$ and a field $F \sbs E$ that is abelian over $Q$ and contains
all roots of unity in $E$. Our criterion (given in the next section) for $L$ to
be a repeated radical extension of $Q$ will be expressed in terms of the
fields $E$ and $F$ and certain associated Galois groups. It will be
convenient to standardize our notation in this situation, and so we will
generally write $G = \gal EQ$, $U = \gal EL$ and $N = \gal EF$. Then
$U \sbs G$ and $N \nor G$ and we define $M = U \cap N$, so that
$M \nor U$. We observe that $|G:U| = |L:Q|$ and that
$|G:UN| = |L \cap F:Q|$. Furthermore, $G/N \cong \gal FQ$ is abelian.

We shall also need to consider finite subgroups of the multiplicative
group $E^\times$ of the field $E$. Note that such subgroups are uniquely
determined by their order. If $D \sbs E^\times$ is a subgroup of order
$n$, for example, then $D$ is exactly the subgroup $\gen\delta$, where
$\delta$ is a primitive $n\,$th root of unity in $E$. Since $G = \gal EQ$
acts on the cyclic group $D$ and $U \sbs G$, we can view $D$
as a $U$-group: a group acted on by $U$.

For the remainder of this section, we focus on the case where $L$
has prime degree $p$ over $Q$. In this situation, of course, $L$ cannot
be a {\it repeated} radical extension of $Q$ unless it is actually a
radical extension. Also, since by the quadratic formula, every degree $2$
extension of fields of characteristic different from $2$ is radical,
we need only consider primes $p > 2$. The following is our principal
result in this situation.
\medskip

\iitem{(2.5) THEOREM.}~~{\sl Let $Q \sbs L$ with $|L:Q| = p$, where
$p$ is an odd prime not equal to the characteristic of $Q$. Let
$E \sps L$ be Galois over $Q$ and suppose that $F \sbs E$ is abelian
over $Q$ and contains the roots of unity in $E$. Let $G$ and its
subgroups, $N$, $U$ and $M$ be as described above. Then $L$ is a
radical extension of $Q$ not contained in $F$ if and only if the following
hold.
\smallskip
\ritem{(i)} $|N:M| = p$.
\smallskip
\ritem{(ii)} $M \nor N$.
\smallskip
\ritem{(iii)} $N/M$ is $U$-isomorphic to a subgroup of $E^\times$.}
\medskip

Given $Q \sbs L$ as in Theorem~2.5, where $|Q:L|$ is prime and not
equal to the characteristic, we see that $L$ is separable over $Q$ and
thus it really is possible to choose $E$ as in the statement of the
theorem. Once $E$ is selected, we have seen that it is easy to find an
appropriate field $F$.

We will apply Theorem~2.5 only in the case where the field $L$ is
quasireal. As we shall see, $L$ cannot be contained in $F$ in that
situation, and this yields  a slight simplification of the result. Before
proceeding with the proof of 2.5, we present the quasireal version as
a corollary. (Recall from Section~1 that a quasireal field is a field of
characteristic zero in which the only roots of unity are $\pm1$.) 
\medskip

\iitem{(2.6) COROLLARY.}~~{\sl In the situation of Theorem~2.5,
assume that $L$ is quasireal. Then $L$ is radical over $Q$ if and only
if conditions (i), (ii) and (iii) hold.}
\medskip

\iitem{Proof.}~~By Theorem~2.5, all that must be proved is that if $L$
is radical over $Q$ and is contained in $F$, then $L$ cannot be quasireal.
But $F$ is abelian over $Q$, and so $L$ is Galois over $Q$ and thus by
Lemma~2.3(a), it is not quasireal.\qed
\medskip

Our proof of Theorem~2.5 relies on a standard fact from group theory,
which may be unfamiliar to some of the readers of this paper. For
completeness, therefore, we state and prove the relevant lemma, and
we mention two of its applications. Recall that a set $\PP$ of subgroups
of a group $G$ is a {\bf partition} of $G$ if $\bigcup\PP = G$ and
$H \cap K = 1$ for distinct members $H,K \in \PP$.
\medskip

\iitem{(2.7) LEMMA.}~~{\sl Let $\PP$ be a partition of a finite group $G$,
and suppose that $G$ acts via automorphisms on an abelian group $A$.
If $A$ contains an element with order not dividing $|\PP| - 1$, then
$\cent AH > 1$ for some member $H \in \PP$.}
\medskip

\iitem{Proof.}~~Write $A$ additively and fix an element $a \in A$ with
order not dividing $|\PP| - 1$. For each subgroup $X \sbs G$, define
$a_X = \sum_{x \in X} a^x$ and note that $a_X \in \cent AX$. We can
assume, therefore, that $a_H = 0$ for all members $H \in \PP$, and also
that $a_G = 0$. Since $\PP$ is a partition of $G$, however, this yields
$$
0 = \sum_{H \in \PP} a_H = (|\PP| - 1)a + a_G = (|\PP| - 1)a \, .
$$
This contradicts our choice of $a$ and completes the proof.\qed
\medskip

One consequence of Lemma~2.7 is the following. If an elementary
abelian $p$-group $P$ of order $p^2$ acts on a nonzero vector space
over a field of characteristic different from $p$, then some subgroup of
order $p$ in $P$ has nontrivial fixed points on the vector space. This is
immediate from the lemma since $P$ is partitioned by its $p+1$
subgroups of order $p$.

Another application of Lemma~2.7 that we shall need concerns a
Frobenius group $F$ with kernel $N$ and complement $H$.
If $F$ acts on a nontrivial vector space in characteristic not dividing
$|N|$, then either $N$ or $H$ must have nontrivial fixed points. This
follows since $F$ is partitioned by $N$ and the $|N|$ conjugates of $H$.
If $N$ has no nontrivial fixed points, then by the lemma, some conjugate
of $H$ has nontrivial fixed points, and it is immediate that $H$ also must
have nontrivial fixed points.
\medskip

\iitem{Proof of Theorem 2.5.}~~Suppose first  that $L \not\sbs F$ and
that $L$ is a radical extension of $Q$. Since $|L:Q|$ is prime, we have
$L \cap F = Q$, and hence $UN = G$.  It follows that
$|N:M| = |G:U| = |L:Q| = p$, proving (i). Furthermore, all roots of unity in
$L$ are in $F$, and hence they lie in $Q$. Since $|L:Q| = p$ and we
are assuming that $L$ is radical over $Q$, we can apply Lemma~2.1(c)
and write $L = Q[\alpha]$, where $\alpha^p \in Q$.

Choose $\tau \in N - M$ and note that since $\tau \not\in U$, we have
$\alpha^\tau \ne \alpha$. But $\tau$ fixes $\alpha^p \in Q$, and it follows
that $\alpha^\tau = \epsilon\alpha$ for some primitive $p\,$th root of unity
$\epsilon \in E$. Write $D = \gen\epsilon \sbs E^\times$, and note that
$|D| = p$.

Since $L[\epsilon]$ is a splitting field over $Q$ for the polynomial
$X^p - \alpha^p \in Q[X]$, we see that $L[\epsilon]$ is Galois over $Q$,
and hence the compositum $L[\epsilon]F$ is also Galois over $Q$.
However, $\epsilon \in F$, and thus $L[\epsilon]F = LF$, and this
corresponds to the subgroup $U \cap N  = M$. It follows that $M \nor G$,
and in particular, $M \nor N$ and (ii) is proved.

Since $|N/M| = p = |D|$, we see that to prove (iii), it suffices to check that
the actions of an arbitrary element $\sigma \in U$ on $D$ and on $N/M$
agree. Recall that we have $\tau \in N - M$ with
$\alpha^\tau = \epsilon\alpha$. The coset $M\tau$ generates $N/M$, and
thus $\tau^\sigma \equiv \tau^s$ mod $M$ for some integer $s$. Also,
$s$ determines the action of $\sigma$ on $N/M$, and since $\epsilon$
generates $D$, it suffices to show that $\epsilon^\sigma = \epsilon^s$.
Now $\epsilon$ is fixed by $\tau$ (because $\epsilon \in F$ and
$\tau \in N = \gal EF$) and $\alpha$ is fixed by $U$, which contains both
$M$ and $\sigma$. We can now compute that
$$
\alpha\epsilon^\sigma = (\alpha\epsilon)^\sigma =
\alpha^{\tau\sigma} = \alpha^{\sigma\inv\tau\sigma} =
\alpha^{\tau^\sigma} = \alpha^{\tau^s} = \alpha\epsilon^s \, .
$$
Thus $\epsilon^\sigma = \epsilon^s$, as desired, and hence $D$ and
$N/M$ are $U$-isomorphic and (iii) holds.

Conversely now, assume the three conditions. By (i), we have
$|UN:U| = |N:M| = p = |L:Q| = |G:U|$, and thus $UN = G$. It follows that
$L \cap F = Q$, and in particular, $L \not\sbs F$. It remains to show that
$L$ is radical over $Q$.

Conditions (ii) and (iii) tell us that $E^\times$ has a subgroup $D$ that is
$U$-isomorphic to $N/M$, which we know has order $p$. In particular,
$D = \gen\epsilon$, where $\epsilon$ is a primitive $p\,$th root of unity
in $E$. Write $K = L[\epsilon]$ and $C = \gal EK$. Thus $C$ is exactly the
set of elements in $U = \gal EL$ that fix $\epsilon$, and hence $C$ is the
kernel of the action of $U$ on $D$, and in particular, $C \nor U$. Also,
$\epsilon \in F$, and thus $M$ fixes $\epsilon$, and we have $M \sbs C$.
By the $U$-isomorphism between $D$ and $N/M$, we see that $C$
acts trivially on $N/M$. Thus $[C,N] \sbs M  \sbs C$, and therefore $N$
normalizes $C$. But $UN = G$ and $U$ normalizes $C$, and we deduce
that $C \nor G$, and hence $K$ is Galois over $Q$.

Now $G$ induces $Q$-linear transformations on $K$, and we write
$\overline G$ to denote the image of $G$ in the full general linear group
$\Gamma = GL_Q(K)$. The map $\sigma \mapsto \overline\sigma$ is a
homomorphism from $G$ onto $\overline G$, and its kernel is
$\gal EK = C$. Since $C \cap N = M$, we see that $\overline N$ has
order $p$ and also that $\overline U$ acts faithfully on $\overline N$ by
conjugation. In particular, $\overline G = \overline U\overline N$ is a
Frobenius group with complement $\overline U$ and kernel
$\overline N$. (Technically, the definition of a ``Frobenius group"
requires that the complement should be nontrivial, which may not be
the case in our situation.)

Since $D \sbs K$, there is another subgroup of order $p$ in $\Gamma$
that is of interest to us, a subgroup different from $\overline N$. This is the
group $\Delta$ consisting of scalar multiplications on $K$ by elements of
$D$. (Note that $\overline N$ fixes the element $1 \in K$ while $\Delta$
does not, and thus $\overline N$ and $\Delta$ really are different.) If
$\mu \in \Delta$ is scalar multiplication by $\delta \in D$ and $a$ is
any element of $K$, then for $\sigma \in G$, we have
$$
(a)\mu^\sigma = (a)\sigma\inv\mu\sigma =
(a^{\sigma\inv}\delta)^\sigma = a\delta^\sigma \, .
$$
It follows that $\overline G$ normalizes $\Delta$ in $\Gamma$ 
and this calculation also shows that $\Delta$ and $D$ are isomorphic
as $G$-groups. Since $D$ and $\overline N$ are $U$-isomorphic by
hypothesis, we conclude that $\Delta$ and $\overline N$ are
$\overline U$-isomorphic.

Recall that $\overline N$ normalizes and is distinct from $\Delta$. It
follows that $\overline N\Delta$ is an elementary abelian subgroup of
$\Gamma$ having order $p^2$, and we write $A = \overline N\Delta$.
Since $\overline U$ acts in the same way on each of $\Delta$ and
$\overline N$, we deduce that every subgroup of $A$ is
$\overline U$-invariant and hence is normal in
$\overline UA = \overline G\Delta$.

Since we are assuming that $Q$ does not have characteristic $p$,
Maschke's theorem applies and we see that $K$ is completely reducible
as a $QN$-module. We can thus write the $Q$-space $K$ as the direct
sum of the subspace consisting of the $N$-fixed points and a unique
complementary $N$-invariant $Q$-subspace $V$, on which $N$ acts
without fixed points. Furthermore, since $N$ acts nontrivially on $K$,
we see that $V$ is nonzero. (Of course, the $Q$-subspace $V \sbs K$ is
not a sub{\it field} of $K$ because $1 \not\in V$.) Since
$\overline N \nor \overline UA$, the uniqueness of $V$ guarantees that
$V$ is invariant under $\overline UA$.

As $A$ is noncyclic of order $p^2$ and acts on $V > 0$, there must be
some subgroup $B \sbs A$ of order $p$ such that $B$ has nontrivial
fixed points on $V$. (This is one of our applications of Lemma~2.7.)
Let $W \sbs V$ be the (nonzero) fixed-point space of $B$. We know that
$B \nor \overline UA$, and it follows that $W$ is invariant under
$\overline UA$. (This is the key point where we use the assumption that
the actions of $U$ on $N/M$ and $D$ agree. It is this assumption that
underlies the fact that every subgroup of $A$ is normalized by
$\overline U$.)

The Frobenius group $\overline U\,\overline N$ acts on $W$, and
$\overline N$ has no nonzero fixed points in $W$ because $W \sbs V$
and $\overline N$ has no nonzero fixed points in $V$. It
follows via Lemma~2.7 that there exists a nonzero element
$\alpha \in W$ fixed by $\overline U$, and thus $\alpha$ is fixed by $U$,
and $\alpha \in L$. Also, since $\alpha \in V$, we know that $\alpha$ is
not fixed by $\overline N$, and thus $\alpha \not\in Q$. Therefore
$L = Q[\alpha]$, and it suffices to show that $\alpha^p \in Q$ in order to
prove that $L$ is a radical extension of $Q$.

Recall that $\alpha \in W$ is fixed  by the subgroup
$B \sbs A = \overline N\Delta$. But $B \ne \Delta$ since $\Delta$ has no
nonzero fixed points on $K$, and it follows that $B$ contains some
element of the form $b = \overline\tau\mu$, where $\overline\tau$
generates $\overline N$ and $\mu \in \Delta$ is multiplication by some
$p\,$th root of unity $\delta$.
We have
$$
\alpha = (\alpha)b = (\alpha)\overline\tau\mu = (\alpha^\tau)\delta \, ,
$$
and thus $\alpha^\tau = \alpha\delta\inv$. We deduce that $\tau$ fixes
$\alpha^p$, which is therefore fixed by all of $\overline N$. Since $U$
fixes $\alpha$, it also fixes $\alpha^p$, and we conclude that $\alpha^p$
is fixed by $\overline U\overline N = \overline G$. Thus $\alpha^p$ is
fixed by $G$, and the proof is complete.\qed
\bigskip

\iitem{3. Characterizing quasireal repeated radical extensions.}

The following is our principal result in this section.
\medskip

\iitem{(3.1) THEOREM.}~~{\sl Suppose $Q \sbs L$, where
$L$ is quasireal. Let $E \sps L$ be Galois over $Q$ and suppose
that $F \sbs E$ is abelian over $Q$ and contains all roots of unity
in $E$. As usual, write $G = \gal EQ$ and let $N = \gal EF$,
$U = \gal EL$ and $M = N \cap U$. Then $L$ is a repeated
radical extension of $Q$ if and only if there is a chain of $U$-invariant
subgroups $M_i$, where $M = M_0 \sbs M_1 \sbs \cdots \sbs M_r = N$,
and all of the following hold.
\smallskip
\ritem{(i)} $|F \cap L:Q|$ is a power of $2$. 
\smallskip
\ritem{(ii)} Each index $|M_i:M_{i-1}|$ is prime.
\smallskip
\ritem{(iii)} $M_{i-1} \nor M_i$ for each integer $i$ with $0 < i \le r$.
\smallskip
\ritem{(iv)} Each factor $M_i/M_{i-1}$ is $U$-isomorphic to a
subgroup of $E^\times$.
\medskip}

In order to understand these conditions better, it seems worthwhile to
play with them a little before we proceed with the proof. Assuming (ii), (iii)
and (iv), let $D_i \sbs E^\times$ be $U$-isomorphic to $M_i/M_{i-1}$. We
see that $|D_i| = |M_i/M_{i-1}|$ is prime, and so the subgroups $D_i$ are
exactly the groups $\gen\delta$ as $\delta$ runs over primitive $p\,$th
roots of unity in $E$ for prime divisors $p$ of $|N:M|$. In particular, (iv)
guarantees the existence in $E$ of all of these $p\,$th roots of unity. To
see exactly which primes these are, observe that
$|N:M| = |NU:U| = |L:L \cap F|$. If (i) holds, then the odd prime divisors
of $|N:M|$ are exactly the odd prime divisors of $|L:Q|$, and this shows
that a consequence of the four conditions is that $E$ contains a
primitive $p\,$th root of unity for each prime divisor $p$ of $|L:Q|$.

Next, observe that $G/N \cong \gal FQ$, which is abelian, by assumption.
Assuming (ii) and (iii), we see that the factor groups $M_i/M_{i-1}$ are
abelian, and thus successive terms of the derived series of $G$ are
contained in the subgroups $M_i$ with decreasing subscripts $i$. In
particular, this tells us that the terms of the derived series of $G$
eventually lie within $M$, and hence $G/H$ is solvable for every normal
subgroup $H$ of $G$ with $H \sps M$. In particular, this holds for all
normal subgroups of $G$ that contain $U$. Translating this last
conclusion into field theory, we see that (ii) and (iii) guarantee that if $Q
\sbs K \sbs L$ and $K$ is Galois over $Q$, then $\gal KQ$ is solvable. Of
course, this is exactly what we would expect by Galois' theorem if $L$
really is a repeated radical extension of $Q$.

We begin working toward a proof of Theorem~3.1 with the following
easy lemma. (A weaker form of this result is Theorem~22.14 of \ref\book,
which appears there with an unnecessarily complicated proof.)
\medskip

\iitem{(3.2) LEMMA.}~~{\sl Suppose that $L$ and $S$ are respectively
a radical extension and a Galois extension of some field $Q$. If both
$L$ and $S$ are quasireal, then $|L \cap S:Q| \le 2$.}
\medskip

\iitem{Proof.}~~Write $L = Q[\alpha]$, where some power of $\alpha$ is
in $Q$, and let $m = |L:L \cap S|$. Since $L$ is quasireal, the only roots
of unity it contains are $\pm1$, and these, of course, lie in $L \cap S$. It
follows by Lemma~2.1(c) that $\alpha^m \in L \cap S$ and we set $\beta
= \alpha^m$ and $F = Q[\beta] \sbs L \cap S$. Observe that $\alpha$ is a
root of the polynomial $f(X) = X^m - \beta \in F[X]$, and thus we have
$$
\deg f = m = |L:L \cap S| = |F[\alpha]:L \cap S| \le |F[\alpha]:F| =
\deg{\minp F\alpha} \le \deg f \, .
$$
Equality must hold throughout, and we deduce that
$L \cap S = F = Q[\beta]$. 

Write $g = \minp Q\beta$ and note that $g$ splits over $S$ since by
hypothesis, $S$ is Galois over $Q$. But some power of $\alpha$ lies in
$Q$, and thus the same is true for $\beta$, and it follows by
Lemma~2.1(a) that every root of $g$ in $S$ has the form $\epsilon\beta$
for some root of unity $\epsilon \in S$. As $S$ is quasireal, the only
possibilities are $\epsilon = \pm1$, and thus $g$ has at most two roots.
The roots of $g$ are distinct, however, and it follows that
$|L \cap S:Q| = \deg g \le 2$.\qed
\medskip

The following result is closely related to Theorems~22.12 and~22.15 of
\ref\book.
\medskip

\iitem{(3.3) THEOREM.}~~{\sl Let $Q \sbs E$, where $E$ is quasireal.
Suppose that $L$ and $S$ subfields of $E$ that are respectively a
repeated radical extension and a Galois extension of $Q$. Then
$L \cap S$ is a repeated quadratic extension of $Q$.}
\medskip

\iitem{Proof.}~~We can assume that $L > Q$, and so we can choose a
radical extension $F$ of $Q$ such that $Q < F \sbs L$. Since $S$ is
Galois over $Q$, it follows by Lemma~3.2 that $|F \cap S:Q| \le 2$.

Now $L$ is a repeated radical extension of $F$ and the compositum
$FS$ is Galois over $F$. As $|L:F| < |L:Q|$, we can
work by induction on $|L:Q|$ and apply the inductive hypothesis with $F$
in place of $Q$ and $FS$ in place of $S$. We deduce that
$D$ is a repeated quadratic extension of $F$, where we have written
$D = L \cap FS$. In other words, there exists a tower of degree
$2$ field extensions $F = F_0 \sbs F_1 \sbs \cdots \sbs F_m = D$. 

Since $S$ is Galois over $Q$, we can apply the so-called ``theorem on
natural irrationalities" (see Theorem~18.22 of \ref\book) to see that
$|FS:F| = |S:F \cap S|$. More generally, if $X$ is any field such that
$F \sbs X \sbs FS$, we have $|FS:X| = |S:X \cap S|$, and
we deduce that $|X:F| = |X \cap S:F \cap S|$. If $X \sbs Y$ are two
consecutive fields in the tower $\{F_i\}$ of degree $2$ extensions
running from $F$ to $D$, it follows that $|Y \cap S:X \cap S| = |Y:X| = 2$.
We conclude that the fields $F_i \cap S$ form a tower of degree $2$
extensions running from $F_0 \cap S = F \cap S$ up to
$F_m \cap S = D \cap S = L \cap S$. Since $|F \cap S:Q|$ is at most $2$,
it follows that $L \cap S$ is a repeated quadratic extension of $Q$, as
required.\qed
\medskip

We mention that in the situation of Theorem~3.3, where we are dealing
with fields of characteristic zero, quadratic extensions are automatically
radical extensions, and thus in the notation of the theorem, $L \cap S$
is a repeated radical extension of $Q$.
\medskip

\iitem{(3.4) COROLLARY.}~~{\sl Let $Q \sbs S \sbs L$, where $L$
is quasireal and $S$ is Galois over $Q$. If $L$ is a repeated radical
extension of $Q$, then $|S:Q|$ is a power of $2$. Conversely, if
$|S:Q|$ is a power of $2$, then at least $S$ is a repeated radical
extension of $Q$.}
\medskip

\iitem{Proof.}~~If $L$ is a repeated radical extension of $Q$, then
we can take $E = L$ in Theorem~3.3, and we deduce that $S$ is a
repeated quadratic extension of $Q$, and hence $|S:Q|$ is a power of
$2$. Conversely, if $|S:Q|$ is a power of $2$, then $\gal SQ$
is a $2$-group, and by elementary group theory and Galois theory, we
see that $S$ is a repeated quadratic extension and hence it is a
repeated radical extension of $Q$.\qed
\medskip

We need one further preliminary result.
\medskip

\iitem{(3.5) LEMMA.}~~{\sl Let $Q \sbs L$, where $L$ is quasireal. Then
$L$ is a repeated radical extension of $Q$ if and only if there is a
tower of fields $Q = L_0 \sbs L_1 \sbs \cdots \sbs L_m = L$ such that
the extensions $L_{i-1} \sbs L_i$ are radical extensions of prime degree
for each integer $i$ with $0 < i \le m$.}
\medskip

\iitem{Proof.}~~Since the sufficiency of the condition is obvious, we
assume that $L > Q$ is a repeated radical extension, and we proceed
to construct the fields $L_i$. Working by induction on $|L:Q|$, we see
that it suffices to construct $L_1 \sbs L$ such that $L_1$ is radical
extension of prime degree over $Q$.

Since $L$ is a proper repeated radical extension of $Q$, we can choose
an element $\alpha \in L$ such that $\alpha \not\in Q$ but
$\alpha^n \in Q$ for some positive integer $n$. Choose $\alpha$ so that
$n$ is as small as possible and observe that this forces $n$ to be a prime
number. Now write $L_1 = Q[\alpha]$, so that $L_1$ is radical over $Q$.
Setting $d = |L_1:Q|$, we have $d \le n$ by Lemma~2.1(b). Since $L$ is
quasireal, however, all  roots  of  unity in $L_1$ lie in $Q$, and thus
$\alpha^d \in Q$ by Lemma~2.1(c). By the minimality of $n$, we deduce
that $d = n$, and $d$ is prime, as required.\qed
\medskip

We can now present the proof of Theorem~3.1, which is our
characterization of quasireal repeated radical extensions. Essentially,
the proof proceeds by repeated application of Theorem~2.5 and
Corollary~2.6.
\medskip

\iitem{Proof of Theorem 3.1.}~~Suppose first that $F \cap L = Q$. In this
situation, we show that $L$ is a repeated radical extension of $Q$ if and
only if there is an appropriate chain of subgroups for which conditions (ii),
(iii) and (iv) of the theorem hold. In view of Lemma~3.5, therefore, we
want to show that the three conditions are equivalent to the existence of
a tower of fields from $Q$ to $L$, where each successive extension is
radical and of prime degree.

By Galois theory, we know that every intermediate field $A$ with
$Q \sbs A \sbs L$ corresponds to a subgroup $W$ with
$U \sbs W \sbs G$ such that $|L:A| = |W:U|$. Also, in our situation,
where $UN = G$, there is a bijective correspondence between
subgroups $W$ with $U \sbs W \sbs G$ and $U$-invariant
subgroups $R$ with $M \sbs R \sbs N$. (Subgroups $W$ and $R$
correspond if $R = W \cap N$, or equivalently $W = UR$.) If $W$ and
$R$ correspond in this situation, we have $|W:U| = |R:M|$.
It follows from all of this that the existence of a chain satisfying
condition (ii) is exactly equivalent to the existence of a tower of fields
from $Q$ to $L$, where each successive extension has prime degree.

Now consider intermediate fields $A$ and $B$ with
$Q \sbs A \sbs B \sbs L$, where $|B:A| = p$, a prime number, and let
$R$ and $S$, respectively, be the corresponding $U$-invariant
subgroups of $N$, so that $M \sbs S \sbs R \sbs N$ and $|R:S| = p$.
Writing $W = \gal EA$ and $V = \gal EB$, we have $R = W \cap N$
and $S = V \cap N = V \cap R$. It suffices to show that $B$
is a radical extension of $A$ if and only if $S \nor R$ and $R/S$ is
$U$-isomorphic to a subgroup of $E^\times$.

If $p = 2$, then $B$ is a quadratic extension of $A$, and this is
automatically a radical extension. Also in this case, $|R:S| = 2$, and so
$S \nor R$ and $R/S$ is $U$-isomorphic to the subgroup
$\gen{-1} \sbs E^\times$. We can thus assume that $p > 2$ and we
appeal to Theorem~2.5 and Corollary~2.6 with $A$ and $B$
in place of the fields $Q$ and $L$ of those results. 

Our present field $E$ is Galois over $A$, and so it will serve as the field
called $E$ in Theorem~2.5. For the field $F$ of 2.5, we take the
compositum $AF$. (This is abelian over $A$ by the theorem on natural
irrationalities, and it certainly contains all roots of unity in $E$.) In
Theorem~2.5, we had $U = \gal EL$, and the corresponding group in the
present situation is $V = \gal EB$. Also in 2.5 we had $N = \gal EF$, and
here, this corresponds to $\gal E{AF} = \gal EA \cap N = R$. Finally, the
group $M$ of 2.5 was $U \cap N$, and the corresponding group here is
$V \cap R = S$. We can thus apply Corollary~2.6 with $V$, $R$ and $S$
in place of $U$, $N$ and $M$, respectively, and it follows that $B$ is
radical over $A$ if and only if $S \nor R$ and $R/S$ is $V$-isomorphic to
a subgroup of $E^\times$. All that remains in this case, therefore, is to
show that $R/S$ is $U$-isomorphic to $\gen\epsilon$ if and only if it is
$V$-isomorphic to $\gen\epsilon$, where $\epsilon$ is a primitive
$p\,$th root of unity in $E$. Since $U \sbs V$, we see that if $R/S$ and
$\gen\epsilon$ are $V$-isomorphic, they are automatically
$U$-isomorphic. To prove the converse, we observe that $S$ acts
trivially on $R/S$ and also that $S$ acts trivially on $\gen\epsilon$ since
$\epsilon \in F$ and $S \sbs N = \gal EF$. But $V = US$, and so we see
that if $R/S$ and $\gen\epsilon$ are $U$-isomorphic, they must also be
$V$-isomorphic, as required.

Finally, we consider the general case, where we do not assume that
$L \cap F = Q$. Applying the previous argument with $L \cap F$ in place
of $Q$ (and $UN$ in place of $G$) we see that $L$ is a repeated
radical extension of $L \cap F$ if and only if conditions (ii), (iii) and (iv)
hold. Observe that $L \cap F$ is a Galois extension of $Q$ because
$F$ is abelian over $Q$. If condition (i) holds, so that $|L \cap F : Q|$
is a power of $2$, then $L \cap F$ is a repeated radical extension of $Q$
by Corollary~3.4. If all four conditions hold, therefore, $L$ is a repeated
radical extension of $L \cap F$, which is a repeated radical extension of
$Q$, and thus $L$ is a repeated radical extension of $Q$, as required.
Conversely,  if $L$ is a repeated radical extension of $Q$, then condition
(i) holds by Corollary~3.4. Also $L$ is a repeated radical extension of
$L \cap F$ in this case, and thus (ii), (iii) and (iv) hold. This completes
the proof.\qed
\bigskip

\iitem{4. Theorem~A.}

Before we proceed to prove Theorems~B, C~and~D, which are our
main applications of Theorem~3.1, we offer an easy proof of the known
Theorem~A. Our proof relies on some of the preliminary
results in Section~3, but it is independent of Theorem~3.1.

The following includes Theorem~A and generalizes it to quasireal fields.
\medskip

\iitem{(4.1) THEOREM.}~~{\sl Let $Q \sbs E$, where $E$ is quasireal,
and suppose $f \in Q[X]$ is irreducible and splits over $E$. If some root
of $f$ lies in a repeated radical extension of $Q$ contained in $E$, then
the splitting field for $f$ over $Q$ in $E$ is a repeated radical extension
of $2$-power degree over $Q$. Also, $\deg f$ is a power of $2$.}
\medskip

\iitem{Proof.}~~Let $L \sbs E$ be a real repeated radical extension of
$Q$ containing a root $\alpha$ of $f$. Then $\alpha \in L \cap S$, where
$S$ is the splitting field for $f$ over $Q$ in $E$. By Theorem~3.3, we
know that $L \cap S$ is a repeated quadratic extension of $Q$, and it
follows that for each element $\sigma \in \gal SQ$, the field
$(L \cap S)^\sigma$ is also a repeated quadratic extension of $Q$. The
compositum $J = \gen{(L \cap S)^\sigma|\sigma \in \gal SQ}$ is
therefore also a repeated quadratic extension of $Q$, and since $J$
contains all the roots $\alpha^\sigma$ of $f$ in $S$, we see that $J = S$
and $S$ is a repeated quadratic extension of $Q$. In particular,
$S$ is a repeated radical extension and the result follows.\qed
\medskip

In the situation of Theorem~4.1, where we have a Galois extension of
$Q$ having $2$-power degree, it is easy to see (without appealing to
Theorem~D and without assuming realness) that every intermediate
field is a repeated radical extension of $Q$. This is because given any
subgroup $H$ of a finite $2$-group $G$, one can always find a chain
of subgroups from $H$ to $G$, each of index $2$ in the next. Every
intermediate field, therefore, is a repeated quadratic extension of $Q$,
and if the characteristic is different from $2$, it is a repeated radical
extension.
\bigskip

\iitem{5. Some group theory.}

We shall need the following result for our proof of Theorem~B.
\medskip

\iitem{(5.1) LEMMA.}~~{\sl Let $M \sbs R \sbs N$ with $M \snor N$,
where $N$ is a finite group. Suppose that $\sigma \in \aut N$ has order
$2$ and that $M$ is $\sigma$-invariant. If $\sigma$ acts
fixed-point-freely on each factor in some $\sigma$-invariant subnormal
series from $M$ up to $N$, then $R \snor N$.}
\medskip

We mention that if a group $H$ acts via automorphisms on a finite group
$N$, and $H$ stabilizes some subnormal subgroup $M$ of $N$, then
there necessarily exists an $H$-invariant subnormal series from $M$
up to $N$. (This is a consequence of Lemma~5.3, below.) In
Lemma~5.1, therefore, it is not necessary to assume the existence
of the $\sigma$-invariant series from $M$ to $N$; the key hypothesis is
that the action of $\sigma$ on each factor in some such series is
fixed-point free. As we shall see in Lemma~5.4, if the action of $\sigma$
on the factors of a $\sigma$-invariant subnormal series from $M$ to
$N$ is fixed-point free, then the same will be true for every such series.

Observe that we did not assume that $R$ is $\sigma$-invariant in
the statement of Lemma~5.1, and in fact, in the situation of that lemma,
one can deduce that $R$ must be $\sigma$-invariant. (When we apply
the lemma in the proof of Theorem~B, however, it will be clear that $R$
is $\sigma$-invariant.)

The hypothesis that $\sigma$ has order $2$ in Lemma~5.1 is
unnecessarily restrictive. If we are willing to assume that $R$ is
$\sigma$-invariant, then the conclusion that $R \snor N$ holds if the
order of $\sigma$ is any prime number. This result is deeper than 5.1,
however, because it relies on J.~Thompson's famous theorem that a
group admitting a fixed-point free automorphism of prime order $p$
must be nilpotent. As is well known, the conclusion of this theorem is
a triviality when $p = 2$, and so the proof of Lemma~5.1 does not
rely on Thompson's result. We have decided, however, to state and
prove the more general theorem.
\medskip

\iitem{(5.2) THEOREM.}~~{\sl Let $\sigma$ be an automorphism of
prime order $p$ of a finite group $N$. Suppose that $M \snor N$ is
$\sigma$-invariant and that $\sigma$ acts fixed-point-freely on each
factor in some $\sigma$-invariant subnormal series from $M$ to $N$.
Let $M \sbs R \sbs N$ and if $p > 2$, assume that $R$ is
$\sigma$-invariant. Then $R \snor N$ and it is $\sigma$-invariant
even when $p = 2$.}
\medskip

The following well-known  result is helpful, but it is not strictly necessary
for the proof of Theorem~5.2. We shall really need this result later,
however.
\medskip

\iitem{(5.3) LEMMA}~~{\sl Let $M \snor N$, where $N$ is a finite group.
Then there exist subgroups $M_i$ such that
$M = M_0 \nor M_1 \nor \cdots \nor M_r = N$ and every automorphism
of $N$ that stabilizes $M$ also stabilizes each of
the subgroups $M_i$.}
\medskip

\iitem{Proof.}~~There is nothing to prove if $M = N$, and so we can
assume that $M \sbs M^N < N$, where $M^N$ denotes the normal
closure of $M$ in $N$. (Note that the normal closure $M^N$ is proper in
$N$ because $M$ is proper and subnormal.) Working by induction
on $|N:M|$, we can find a subgroup chain
$M = M_0 \nor M_1 \nor \cdots \nor  M_{r-1} = M^N$,
where each subgroup $M_i$ is stabilized by the automorphisms of
$M^N$ that stabilize $M$. Since the automorphisms of $N$ that
stabilize $M$ also stabilize $M^N$, the result follows by defining
$M_r = N$.\qed
\medskip

\iitem{(5.4) LEMMA.}~~{\sl Let $M \snor N$ and suppose $M$ is
$\sigma$-invariant, where $\sigma \in \aut N$ has prime order $p$. If
$\sigma$ acts fixed-point-freely on all factors in some $\sigma$-invariant
subnormal series from $M$ up to $N$, then $\sigma$ acts
fixed-point-freely on every section $R/S$ with $M \sbs S \nor R \sbs N$,
where $R$ and $S$ are $\sigma$-invariant.}
\medskip

\iitem{Proof.}~~Let $\XX$ be a $\sigma$-invariant subnormal series
from $M$ to $N$ such that $\sigma$ acts fixed-point-freely on each
factor. We claim that the only $\sigma$-invariant right coset of $M$ in
$N$ is $M$ itself, and thus $|N:M| \equiv 1$ mod $p$. To see this,
suppose that the coset $My$ is $\sigma$-invariant and consider the
minimal term $Y$ of $\XX$ that contains $y$. If $Y = M$, then $My = M$
as desired, and so we suppose that $Y > M$ and derive a contradiction.
Consider the term $X$ just below $Y$ in $\XX$. Then $X \nor Y$, both
$X$ and $Y$ are $\sigma$-invariant and by hypothesis, the action of
$\sigma$ on $Y/X$ is fixed-point free. But $Xy = X(My)$ is
$\sigma$-invariant, and thus $y \in X$, contradicting our choice of $Y$.

Now let $R$ and $S$ be as in the statement of the lemma
and suppose that $\sigma$ stabilizes the coset $Sr \in R/S$. As
$M \sbs S$, we see that $Sr$ is a union of $|S:M|$ right cosets of
$M$, and these are permuted by $\sigma$. But $|S:M|$ divides
$|N:M|$, which is not divisible by $p$. It follows that $\sigma$ fixes
one of the cosets of $M$ in $Sr$, and thus $M \sbs Sr$. We conclude
that $Sr = S$, as required.\qed
\medskip

Under the hypotheses of Theorem~5.2, we see by Lemma~5.4 that
if $M \sbs D \nor N$ for some $\sigma$-invariant subgroup $D$, then
$\sigma$ acts fixed-point-freely on $N/D$, which is therefore nilpotent by
Thompson's theorem. Also, if $p = 2$, then $\sigma$ inverts all elements
of $N/D$, which is therefore abelian.  In this case, every subgroup $R$
satisfying $D \sbs R \sbs N$ is $\sigma$-invariant.

If $M \nor N$ in Theorem~5.2, therefore, then $N/M$ is nilpotent, and
it is immediate that $R \snor N$. The significance of Theorem~5.2 is that
it is not necessary to assume that $M$ is normal; subnormality is
sufficient.
\medskip

\iitem{Proof of Theorem~5.2.}~~The result is trivially true when $R = N$,
and so we can assume that $R < N$, and we work by double induction:
first on $|N|$ and then on $|N:R|$. If there exists a subgroup $S$ with
$R < S < N$, where $S$ is $\sigma$-invariant if $p > 2$, then by the
inductive hypothesis applied to the situation $S \sbs N$, we deduce that
$S \snor N$ and that $S$ is $\sigma$-invariant (even if $p = 2$). By
Lemma~5.3 (or by intersecting the given series with $S$), we see that
there is a $\sigma$-invariant subnormal series from $M$ to $S$. Also, by
Lemma~5.4, the hypotheses apply with $S$ in place of $N$ (with the
same subgroups $M$ and $R$). It follows by the inductive hypothesis
applied in the situation $R \sbs S$ that $R$ is $\sigma$-invariant and is
subnormal in $S$, and we are done in this case. We can thus suppose
that $R$ is a maximal subgroup of $N$ if $p = 2$, and that it is a maximal
$\sigma$-invariant subgroup if $p > 2$.

Let $H = M^N$, the normal closure, and write $D = R \cap H$.
Observe that $H < N$ because $M$ is proper and subnormal in $N$,
and $D \nor R$ since $H \nor N$. Also, $H$ is $\sigma$-invariant
because $M$ is, and $D$ is $\sigma$-invariant if $p > 2$.
By Lemmas~5.3 and~5.4, therefore, we can apply the inductive
hypothesis to the situation $M \sbs D \sbs H$, and we deduce that
$D \snor H$ and that $D$ is unconditionally $\sigma$-invariant. 
It follows that either $D = H$, in which case $D \nor N$, or else
$\norm HD > D$. In the latter situation, $\norm ND \not\sbs R$, and so
$\norm ND > R$. Since $\norm ND$ is $\sigma$-invariant, it follows
from the maximality of $R$ that $\norm ND = N$. In either case,
therefore, we have $D \nor N$.

By Lemma~5.4, the action of $\sigma$ on $N/D$ is fixed-point free,
and thus $N/D$ is nilpotent and $R \snor N$. Also, if $p = 2$, then
$N/D$ abelian, and each of its elements is inverted by $\sigma$.
It follows in this case that $R$ is $\sigma$-invariant.\qed
\bigskip

\iitem{6. Theorem B.}

We are finally ready to prove Theorem~B, which we restate here.
\medskip

\iitem{(6.1) THEOREM.}~~{\sl Suppose that $Q$ is a real field and that
$Q \sbs L$ is a repeated radical extension with $|L:Q|$ odd. If
$Q \sbs K \sbs L$, then $K$ is a repeated radical extension of $Q$.}
\medskip

\iitem{Proof.}~~We can assume that $L \sbs \C$, the complex numbers.
Choose a Galois extension $E \sps Q$ with $L \sbs E \sbs \C$, and write
$G = \gal EQ$. Since $Q$ is real, $E$ is invariant under complex
conjugation and we let $\sigma \in G$ be the restriction of conjugation to
$E$. Write $U = \gal EL \sbs G$ and note that $|G:U| = |L:Q|$ is odd, and
thus by Sylow's theorem, some conjugate $U^\tau$ of $U$ in $G$
contains $\sigma$. We can replace $L$ by the $Q$-isomorphic field
$L^\tau$, and we can thus assume that $\sigma \in U$. (Note that the
property of being a repeated radical extension of $Q$ is preserved by
$Q$-isomorphism.) Since $\sigma \in U$, we have $L \sbs \R$, and in
particular, $L$ is quasireal and Theorem~3.1 applies.

Let $F \sbs E$ be abelian over $Q$ and contain all roots of unity in $E$,
and write $N = \gal EF$ and $M = U \cap N$, as usual, so that there is an
appropriate chain of subgroups for which the four conditions of
Theorem~3.1 hold. We will show that $K$ is a repeated
radical extension of $Q$ by verifying these conditions for $K$. We
thus define $V = \gal EK \sps U$ and $R = V \cap N \sps M$, and we
work with $V$ and $R$ in place of $U$ and $M$.

By (i), we know that $|L \cap F:Q|$ is a power of $2$. But $|L:Q|$ is
odd, by hypothesis, and thus $L \cap F = Q$ and $K \cap F = Q$, and
so (i) holds for the field $K$. Also $UN = G$ in this situation, and thus
$|N:M| = |G:U| = |L:Q|$ is odd. We proceed to verify (ii), (iii) and (iv) for
$K$.

Conditions (ii), (iii) and (iv) for $L$ tell us that $M \snor N$ and that there
is a $U$-composition series $\XX$ for $N$ through $M$ such each factor
of $\XX$ above $M$ is $U$-isomorphic to a group of roots of unity of
prime order. As $|N:M|$ is odd, these primes are all odd, and thus
complex conjugation acts fixed-point-freely on each of these groups of
roots of unity. Since $\sigma \in U$, it follows that $\sigma$ also acts
fixed-point-freely on the factors of $\XX$ above $M$, and Lemma~5.1
applies. We conclude that $R \snor N$.

Now $R \nor V$, and so by Lemma~5.3, we can construct a
$V$-invariant subnormal series from $R$ to $N$, and this can be refined
to a $V$-composition series $\YY$ for $N$ that has $R$ as one of its
terms. Observe that $V = UR$ and, of course, $R$ acts trivially on each
the factors of $\YY$ above $R$. These factors are therefore $U$-simple,
and hence they are $U$-isomorphic to some of the factors above $M$ in
the $U$-composition series $\XX$. In particular, each factor $Y$ of
$\YY$ above $R$ is $U$-isomorphic to some subgroup
$D \sbs E^\times$ of prime order. Conditions (ii) and (iii) of
Theorem~3.1thus hold.

To complete the proof, it suffices to show that $Y$ and $D$ are
actually $V$-isomorphic. But $D$ consists of roots of unity, and so
$D \sbs F$ and $N$ acts trivially on $D$. In particular, $R$ acts trivially
on $D$. Since $R$ also acts trivially on $Y$ and $V = UR$, it follows
that $Y$ and $D$ are $V$-isomorphic, as required.\qed
\bigskip

\iitem{7. Theorem C.}

Let $f \in Q[X]$ be irreducible, where $Q$ is a real field and $f$ has at
least one root that is contained in a real repeated radical extension of
$Q$. By Theorem~A, we know that only when $\deg f$ is a power of $2$
can it be true that all of the complex roots of $f$ are real.  In the opposite
extreme case, where $\deg f$ is odd, Theorem~C asserts that $f$ can
have only the one real root with which we started.  We are now ready to
prove this.
\medskip

\iitem{Proof of Theorem~C.}~~We are given that $f$ has a root $\alpha$
lying in some real repeated radical extension $L$ of $Q$, and we choose
a field $E$, Galois over $Q$, with $L \sbs E \sbs \C$. Let $F \sbs E$ be
abelian over $Q$ and contain all roots of unity in $E$, and let $G$, $U$,
$N$ and $M$ have their usual meanings, so that the four conditions of
Theorem~3.1 hold for an appropriate chain of subgroups. Note that $f$
splits over $E$, and our task is thus to show that $\alpha$ is the only
real root of $f$ in $E$.

We argue first that it is no loss to assume that $L \cap F = Q$.
To see why this is so, write $K = L \cap F$ and note that
$Q[\alpha] \sbs K[\alpha]$. It follows that $\deg f = |Q[\alpha]:Q|$
divides $|K[\alpha]:Q| = |K[\alpha]:K||K: Q|$. Since $\deg f$ is odd
and $|K:Q|$ is a power of $2$ by the first condition of Theorem~3.1, we
deduce that $\deg f$ divides $|K[\alpha]:K|$. It follows from this that $f$
is irreducible over $K$. Since $L$ is a real repeated radical extension of
$K$, we can replace our ground field $Q$ with $K$, leaving $E$ and $F$
unchanged. (Note that $E$ is Galois over $K$ and that $F$ is abelian
over $K$.) We can thus assume that $L \cap F = Q$, as claimed, and so
we have $UN = G$.

Let $\beta \in E$ be a real root of $f$ and recall that we must show
that $\beta = \alpha$. Since $G = UN$ acts transitively on the roots of $f$
in $E$ and $U$ fixes $\alpha$ (because $\alpha \in L$), there exists an
element of $N$ that carries $\alpha$ to $\beta$.

By Theorem~3.1, we have a $U$-invariant subnormal series $\XX$
from $M$ to $N$ with factors $U$-isomorphic to prime-order subgroups
of $E^\times$. Let $X$ be the least term in $\XX$ that contains an
element $\tau$ carrying $\alpha$ to $\beta$. If $X = M$, then
$\tau \in M \sbs U = \gal EL$, and $\tau$ fixes $\alpha \in L$. In this case,
$\beta = \alpha^\tau = \alpha$, as required. We can thus assume that
$X > M$, and we let $Y$ be the term just below $X$ in the series $\XX$.
In particular, $M \sbs Y \nor X$ and $X/Y$ is $U$-isomorphic to
$\gen\epsilon$, where $\epsilon$ is a primitive $p\,$th root of unity in
$E$ for some prime $p$.

Now let $\sigma \in G$ be the restriction of complex conjugation to $E$
and note that $\sigma \in U$ since $L$ is real. Thus $X^\sigma = X$ and
in particular, $\tau^\sigma \in X$ and $\tau^\sigma\tau\inv \in X$. Also,
since $\beta$ and $\alpha$ are both real, we compute that
$$
\alpha^{\sigma\tau\sigma} = \alpha^{\tau\sigma} = \beta^\sigma =
\beta = \alpha^\tau \,.
$$
Thus $\tau^\sigma\tau\inv$ fixes $\alpha$, and hence it lies in
$X_\alpha$, the stabilizer in $X$ of $\alpha$.

By the minimality of $X$, we see that no element of $Y$ carries $\alpha$
to $\beta$, and thus we cannot have $X = X_\alpha Y$. Since $Y \nor X$
has prime index, we deduce that $X_\alpha \sbs Y$. We know,
however, that $|G:G_\alpha| = \deg f$ is odd, and since $X \snor G$, it
follows that $|X:X_\alpha|$ is odd, and thus $|X/Y|$ is odd.

Since $\sigma \in U$ inverts the elements of $\gen\epsilon$, we deduce
that $\sigma$ inverts the elements of $X/Y$, and therefore no nonidentity
element of $X/Y$ is fixed by $\sigma$. But
$\tau^\sigma \tau\inv \in X_\alpha \sbs Y$, and this shows that the coset
$Y\tau$ is a $\sigma$-fixed point of $X/Y$. We conclude that $\tau \in Y$,
and this contradicts the choice of $X$.\qed
\bigskip

\iitem{8. Theorem D.}

As another application of Theorem~3.1 we deduce Theorem~D. In fact,
as promised, we have the following slightly stronger result.
\medskip

\iitem{(8.1) THEOREM.}~~{\sl Let $Q \sbs L$ be a repeated radical
extension, where $L$ is quasireal. If $|L:Q|$ is a prime power and
$Q \sbs K \sbs L$, then $K$ is a repeated radical extension of $Q$.}
\medskip

Because Theorem~B deals with the case where $|L:Q|$ is odd, our
primary interest in Theorem~8.1 is in the situation where $|L:Q|$ is a
power of $2$, as in Theorem~D. The cases of 8.1 where $|L:Q|$ is an
odd prime power are not completely covered by Theorem~B, however,
because we assume here only that $L$ is quasireal, while in
Theorem~B, the assumption is that $Q$ is actually real.

We need the following easy lemma from group representation theory.
\medskip

\iitem{(8.2) LEMMA.}~~{\sl Let $V$ be a finite group and suppose that
$U \sbs V$ is a subgroup, where $|V:U|$ is a power of a prime number
$p$. Let $X$ be a simple $FV$-module, where $F$ has characteristic
$p$, and suppose that all composition factors of $X$ viewed as an
$FU$-module are isomorphic and of dimension $1$. Then
$\dimm FX = 1$.}
\medskip

\iitem{Proof.}~~Let $R$ be a Sylow $r$-subgroup of $U$, where
$r \ne p$. Then $X$ is semisimple as an $FR$-module, and since
$R \sbs U$, all composition factors of this $FR$-module have
dimension $1$ and are isomorphic. It follows that each element of
$R$ acts via scalar multiplication on $X$.

Now let $Z \sbs V$ be the subgroup consisting of all elements that
act via scalar multiplication. We have seen that $Z$ contains a full
Sylow $r$-subgroup of $U$ for each prime $r \ne p$, and since 
$|V:U|$ is a power of $p$, it follows that $|V:Z|$ is a power of $p$.
Thus $V = PZ$, where $P$ is some Sylow $p$-subgroup of $V$. 

Now $P$ fixes some nonzero element $x \in X$, and thus both $P$
and $Z$ stabilize the subspace $Fx \sbs X$. Since $PZ = V$ and
$X$ is simple as an $FV$-module, we deduce that $Fx = X$ and the
proof is complete.\qed
\medskip

In the following, we use the standard group-theoretic notation $\Oh pN$
for a finite group $N$. Recall that this is the unique smallest normal
subgroup of $N$ having $p$-power index, where $p$ is a prime number.
It is clear that $\Oh pN \sbs M$ whenever $M \nor N$ and $|N:M|$ is
a power of $p$. In fact, an easy inductive argument shows that
$\Oh pN \sbs M$ whenever $M$ is subnormal in $N$ and has
$p$-power index.
\medskip

\iitem{Proof of Theorem 8.1.}~~To apply Theorem~3.1, let $E \sps L$ be
Galois over $Q$ and suppose that $F \sbs E$ is abelian over $Q$ and
contains all roots of unity in $E$. As usual, write $G = \gal EQ$,
$U = \gal EL$, $N = \gal EF$ and $M = U \cap N$ and note that the
four conditions of Theorem~3.1 must hold since $L$ is a repeated radical
extension of $Q$. To show that $K$ is a repeated radical extension of
$Q$, we let $V = \gal EK \sps U$ and we write $R = V \cap N \sps M$.
We must verify the four conditions of 3.1 in the situation where $K$
replaces $L$, so that $V$ replaces $U$ and $R$ replaces $M$.

Since $|K \cap F:Q|$ divides $|L \cap F:Q|$, which we are assuming is a
power of $2$, the first condition is satisfied for $K$, and we work toward
proving (ii), (iii) and (iv).

Note that $|N:M| = |UN:U|$ divides $|G:U| = |L:Q|$, which is a
power of some prime $p$. Since (iii) holds for the field $L$, we know
that $M$ is subnormal in $N$ with $p$-power index, and it follows that
$\Oh pN \sbs M \sbs R \sbs N$, and thus $R$ is subnormal in $N$.
Also $R$ is $V$-invariant, and thus by Lemma~5.3, there exists a
$V$-composition series $\XX$ for $N$, having $R$ as one of its terms.
Let $X$ be any one of the factors of $\XX$ above $R$. To prove the
three conditions we show that $|X| = p$ and that $X$ is $V$-isomorphic
to the subgroup $\gen\epsilon \sbs E^\times$, where $\epsilon$
is a primitive $p\,$th root of unity in $E$.

Since $U \sbs V$, we see that $R$ is $U$-invariant, and the
$U$-composition factors of $N$ above $R$ are among the
$U$-composition factors of $N$ above $M$. Each of these, however,
is $U$-isomorphic to $\gen\epsilon$, and it follows that when $X$ is
viewed as a $U$-group, all of its composition factors are isomorphic
and of order $p$. Since $|V:U|$ is a power of $p$, we can apply
Lemma~8.2 to deduce that $X$ has order $p$, and thus $X$ is
$U$-isomorphic to $\gen\epsilon$. But each $p$-element of $V$ acts
trivially on both $X$ and $\gen\epsilon$, and since $V$ is generated by
$U$ and $p$-elements, it follows that $X$ and $\gen\epsilon$ are
actually $V$-isomorphic, as required.\qed
\bigskip

\iitem{9. Examples and further remarks.}

Let $Q$ be a real field and suppose that $f \in Q[X]$ is irreducible of
degree $n$ and that $f$ has a root that lies in a real repeated radical
extension of $Q$. By Theorem~A, we know that if $f$ has $n$ real roots,
then $n$ must be a power of $2$. To see that this actually can happen
when $n$ is an arbitrary  power of $2$, let $p$ be any prime congruent to
$1$ modulo $2n$ and let $E$ be the unique extension of degree $n$
over $\Q$ contained in the cyclotomic field of $p\,$th roots of unity. Then
$E$ is a real field and $\gal E\Q$ is a (cyclic) $2$-group. It follows that
$E$ is a repeated radical extension of $\Q$, and so if we take $f$ to be
the minimal polynomial over $\Q$ of any generating element of $E$, we
have the desired example.

If $n = \deg f$ is odd, on the other hand, then Theorem~C tells us that $f$
has only one real root. This suggests that perhaps in general, when $n$
is not necessarily either odd or a power of $2$, the number of real roots
of $f$ is at most the $2$-part of of $n$. This is incorrect, however, and we
give an explicit example of an irreducible polynomial $f \in \Q[X]$ of
degree $6$ having four real roots, of which exactly one lies in a real
repeated radical extension of $\Q$. (This also shows that not all real
roots of an irreducible polynomial over $\Q$ need be ``alike": it is
possible for some to lie in real repeated radical extensions while others
do not.)

We claim that the polynomial $f(X) = (X^3 - 3X + 3)^2 - 3$ has the desired
properties. First, note that $f$ is irreducible over $\Q[X]$ by the
Eisenstein criterion since the constant term of $f$ is $6$, the leading
coefficient is $1$ and all of the other coefficients are divisible by $3$.
Next, we factor $f(X) = (X^3 - 3X + 3 + \sqrt3)(X^3 - 3X + 3 - \sqrt3)$, and
we investigate the (complex) roots of each factor.

Let $x$ be a real variable and consider the polynomial function
$h(x) = x^3 - 3x + a$, where $a$ is a real number. Since $h$ has a local
maximum at $x = -1$ and a local minimum at $x = 1$, we see that the
graph of $y = h(x)$ meets the $x$-axis as many as three times if and only
if $h(-1) > 0$ and $h(1) < 0$. Since $h(-1) = a + 2$ and $h(1) = a - 2$, it
follows that the condition for $h$ to have three real zeros is that
$-2 < a < 2$. (This can also be checked by considering the discriminant 
of $h$.) The number $a = 3 - \sqrt3$ clearly satisfies this condition,
but $a = 3 + \sqrt3$ does not. Returning now to the factors of the
polynomial $f(X)$, we deduce that $u(X) = X^3 - 3X + 3 + \sqrt3$ has
exactly one real root while $v(X) = X^3 - 3X + 3 - \sqrt3$ has three real
roots. Therefore, $f(X)$ has a total of four real roots, as claimed.

Next, we observe that each of the polynomials $u$ and $v$ is irreducible
over $\Q[\sqrt3]$ since otherwise, one of these polynomials, and
therefore also $f$, would have a root in this quadratic extension of $\Q$,
and this is impossible since $f$ is irreducible of degree $6$ over $\Q$. By
Theorem~A, therefore, none of the three real roots of $v$ lies in a real
repeated radical extension of $\Q[\sqrt3]$, and thus none lies in a real
repeated radical extension of $\Q$.

What remains is to show that the unique real root of $u$ does lie in a real
repeated radical extension of $\Q$, and for this purpose, it suffices to
show that it lies in a real repeated radical extension of $\Q[\sqrt3]$. The
following result does the job.
\medskip

\iitem{(9.1) THEOREM.}~~{\sl Let $Q$ be a real field and suppose that
$f \in Q[X]$ is an irreducible cubic polynomial having exactly one real root
$\alpha$. Then $\alpha$ lies in a real repeated radical extension of $Q$.}
\medskip

Note that the converse of Theorem~9.1 is also true: if the irreducible
cubic polynomial $f$ has a root $\alpha$ that lies in a real repeated
radical extension of $Q$, then $\alpha$ is the only real root of $f$. This
case of Theorem~C also follows by Theorem~A, since if a real cubic
polynomial has two real roots, it has three.

Actually, something slightly more general than Theorem~9.1 is true,
and so we state this improved result and prove it instead.
\medskip

\iitem{(9.2) THEOREM.}~~{\sl Let $Q$ be a real field and suppose that
$f \in Q[X]$ is a solvable irreducible polynomial of degree $p$ over $Q$,
where $p$ is a Fermat prime. If $f$ does not split over $\R$, then $f$
has a root that lies in a real repeated radical extension of $Q$.}
\medskip

\iitem{Proof.}~~Let $S$ be the splitting field for $f$ over $Q$ in $\C$ and
let $H = \gal SQ$, so that $H$ is a solvable permutation group of prime
degree $p$. It is well known and easy to prove that $H$ must have a
normal subgroup $P$ of order $p$ and that $H/P$ isomorphic to a
subgroup of the abelian group $\aut P$, which in our case, where $p$ is
Fermat, has $2$-power order. It follows that there exists a field $T$ with
$Q \sbs T \sbs S$, such that $T$ is abelian over $Q$, and where
$|S:T| = p$ and $|T:Q|$ is a power of $2$.

Next, we define $E = S[\epsilon]$, where $\epsilon$ is a complex
primitive $p\,$th root of unity. (It is possible, of course, that $\epsilon \in
S$, in which case $E = S$.) Note that $Q[\epsilon]$ is Galois over $Q$ of
degree dividing $p - 1$, which is a power of $2$. Since
$E = SQ[\epsilon]$ is a compositum of Galois extensions of $Q$,
we see that $E$ is Galois over $Q$. Also, by the natural irrationalities
theorem, $|E:S|$ divides $|Q[\epsilon]:Q|$ and we see that
$|E:Q| = 2^ep$, for some integer $e$.

The restriction of complex conjugation to $E$ is an element $\sigma$
of $G = \gal EQ$, and we write $U = \gen\sigma$ and $L = E \cap \R$,
so that $\gal EL = U$. Since $\deg f = p$ is odd, $f$ has some real root
$\alpha$, and we have $\alpha \in L$. We are assuming that $f$ does not
split over $\R$, and in particular, it does not split over $L$ and $L$ is not
Galois over $Q$. It follows that $U$ is not normal in $G$, and thus
$U$ is nontrivial, so that $|U| = 2$. 

We propose to complete the proof by appealing to Theorem~3.1 to show
that $L$ is a repeated radical extension of $Q$. For this purpose, we
need a field $F \sbs E$ that is abelian over $Q$ and contains all roots of
unity in $E$. Because $T$ is abelian over $Q$, it
follows by Lemma~2.4 that we can choose $F$ so that it contains $T$.
(Actually, it is not hard to see that we can take $F = T[\epsilon]$,
but we shall not need that fact.) Now $\gal FQ$ is abelian, and thus
if $\alpha \in F$, the polynomial $f$ would split over the real field
$Q[\alpha]$. This is contrary to the hypothesis, and therefore,
$\alpha \not\in F$. Thus $T \sbs F \cap S < S$, and since $|S:T| = p$
is prime, we deduce that $F \cap S = T$. Also, $\epsilon \in F$, and thus
$SF = E$, and we conclude by the natural irrationalities theorem
that $|E:F| = |S:T| = p$. Writing $N = \gal EF$, as usual, we see
that $|N| = p$.

We are now ready to check the four conditions of Theorem~3.1. Since
$N \nor G$ has order $p$ and $|U| = 2$, we see that
$|UN| = 2p$.  Thus $|F \cap L:Q| = |G:UN|$, and this is a power of $2$
because $|G| = 2^ep$. This verifies the first condition of Theorem~3.1.

Next, observe that $M = U \cap N$ is trivial, and since $|N| = p$, the
second and third conditions hold for the subgroup chain $M \sbs N$. To
verify the fourth condition, we need to show that $N$ is $U$-isomorphic
to $\gen\epsilon$. The unique nonidentity element $\sigma$ of $U$
inverts the elements of $\gen\epsilon$, and so it suffices to show that
$U$ acts nontrivially on $N$. (The only possible nontrivial action of $U$
on the group $N$ of prime order is for the involution in $U$ to invert all
elements of $N$.) In other words, it is enough to establish that $U$ is
not normal in $UN$. Observe that $UN \nor G$ since
$G/N \cong \gal FQ$ is abelian. If $U \nor UN$, then $U$ would be
characteristic in $UN$, and hence $U \nor G$. We know that this is not
the case, however, and this completes the proof.\qed
\medskip

This completes the argument that the sixth degree irreducible polynomial
over $\Q$ that we described previously does indeed have exactly four
real roots and that exactly one of them is in a real repeated radical
extension of $\Q$.

An obvious question at this point is whether or not the hypothesis
that $p$ is a Fermat prime is really necessary in Theorem~9.2. Although
we have not found an explicit example, it seems likely that the conclusion
of 9.2 does not hold more generally.  

In the situation of Theorem~9.1, we know that $\alpha$ lies in a real
repeated radical extension of $Q$, but it is not necessarily the case that
the cubic extension $Q \sbs Q[\alpha]$ is radical. (And when it is not
radical, it is obviously not a repeated radical extension either.) This
shows that in general, subfields of real repeated radical extensions of
a real field $Q$ need not be repeated radical extensions of $Q$. In
fact, an example exists over the rational numbers $\Q$.
\medskip

\iitem{(9.3) EXAMPLE.}~~{\sl There exists an irreducible cubic
polynomial over $\Q$ having a unique real root $\alpha$, where
$\Q[\alpha]$ is not a repeated radical extension of $\Q$. In fact, the
polynomial $X^3 - 3X + 3$ has this property.}
\medskip

The construction for Example~9.3 relies on the following easy result.
\medskip

\iitem{(9.4) THEOREM.}~~{Let $f \in Q[X]$ be an irreducible cubic
polynomial, where $Q$ is a quasireal field.  Suppose that $\alpha$ is a
root of $f$ in some extension field of $Q$ and that $Q[\alpha]$ is a
repeated radical extension of $Q$. Then the discriminant of $f$ is
$-3m^2$ for some element $m \in Q$.}
\medskip

\iitem{Proof}~~Let $S$ be a splitting field for $f$ over $Q$ and let
$\Delta$ be the discriminant of $f$.  Since $|Q[\alpha]:Q| = 3$,
we see that $Q[\alpha]$ must actually be a radical extension of $Q$, and
so by Lemma~2.1(c), we know that $Q[\alpha] = Q[\beta]$, where
$\beta^3 \in Q$.  The polynomial $X^3 - \beta^3$ is thus irreducible in
$Q[X]$, and hence it splits over $S$.  It follows that $S$ contains the
primitive cube root of unity $\omega = (-1 + \sqrt{-3})/2$, and so in
particular, it contains $\sqrt{-3}$, and $-3$ is a square in $S$. Also,
since $Q$ is quasireal, $\omega \not\in Q$, and thus $\gal SQ$
has order $6$ and is isomorphic to the full symmetric group of degree
$3$. It follows that $\Delta$ is not a square in $Q$. In this case, $S$
contains a unique quadratic extension $T$ of $Q$ and $T$ is
Galois over $Q$. Each of $\sqrt\Delta$ and $\sqrt{-3}$ must lie in $T$
and each of these elements is negated by the unique nonidentity
automorphism in $\gal TQ$. It follows that $\sqrt\Delta/\sqrt{-3}$ is
fixed by $\gal TQ$, and it hence lies in $Q$. In other words, $\Delta/(-3)$
is a square in $Q$, as desired.\qed
\medskip

\iitem{Proof of Example 9.3.}~~Let $f(X) = X^3 - 3X + 3 \in \Q[X]$,
so that $f$ is irreducible by the Eisenstein criterion. Also, since the
constant term of $f$ does not lie between $-2$ and $2$, we know by
our earlier analysis that $f$ must have exactly one real root $\alpha$.

Recall that the the discriminant $\Delta$ of the polynomial
$X^3 + bX + c$ is $-4b^3 - 27c^2$. (See the discussion following
Lemma~23.21 in \ref\book, for example.) In our case, where
$b = -3$ and $c = 3$, we compute that $\Delta(f) = -(5)(27)$, and this
is not of the form $-3m^2$ with $m \in \Q$. It follows by Theorem~9.4
that $Q[\alpha]$ cannot be a repeated radical extension of $\Q$.\qed
\medskip

There is an analog of Theorem~9.1 for quartic polynomials. Although
this result too can be proved using Theorem~3.1, we have decided to
use an alternative approach.
\medskip

\iitem{(9.5) THEOREM.}~~{\sl Let $Q$ be a real field and suppose that
$f \in Q[X]$ is an irreducible quartic polynomial having exactly two real
roots. Then each real root of $f$ lies in a real repeated radical extension
of $Q$.}
\medskip

\iitem{Proof.}~~Let $\alpha$, $\beta$, $\gamma$ and $\delta$ be the
four (distinct) complex roots of $f$, where $\alpha$ and $\beta$ are real
and $\gamma$ and $\delta$ are nonreal complex conjugates.
Define the complex numbers $r = \alpha\beta + \gamma\delta$, 
$s = \alpha\gamma + \beta\delta$ and $t = \alpha\delta + \beta\gamma$,
and observe that $r$ is real and that $s$ and $t$ are distinct since
$s - t = (\alpha - \beta)(\gamma - \delta) \ne 0$. Also, $s$ and $t$ are
complex conjugates, and so they are nonreal.

We claim that $r$ is contained in some real repeated radical
extension $L$ of $Q$. To see why this is so, observe that the Galois
group of $f$ over $Q$ permutes the set $\{r,s,t\}$. This group thus fixes
the coefficients of the polynomial $g(X) = (X - r)(X - s)(X - t)$, and we
deduce that $g \in Q[X]$. If $g$ is reducible over $Q$, then
$|Q[r]:Q| \le 2$, and in this case, $Q[r]$ is a radical extension of $Q$ and
we can take $L = Q[r]$. Otherwise, $g$ is irreducible over $Q$, and
since $r$ is the unique real root of $g$, it follows that $L$ exists by
Theorem~9.1.

Let $u = \alpha\beta\gamma\delta$ and observe that $u \in Q \sbs L$.
We compute that $r\alpha\beta = (\alpha\beta)^2 + u$, and thus
$\alpha\beta$ satisfies a quadratic equation over $L$. Thus
$L[\alpha\beta]$ is a real field of degree at most $2$ over $L$, and
hence it is a real repeated radical extension of $Q$. Replacing $L$
by $L[\alpha\beta]$, therefore, we can assume that $\alpha\beta \in L$.

Now let $v = \alpha+\beta+\gamma+\delta \in Q \sbs L$. Then
$$
(\alpha + \beta)(v - (\alpha + \beta)) =
(\alpha + \beta)(\gamma + \delta) = s + t = (r + s + t) - r \in L \,,
$$
since $r + s + t \in Q$. Thus $\alpha + \beta$ satisfies a quadratic
equation over $L$, and reasoning as before, we can replace $L$ by
$L[\alpha + \beta]$ and assume that $\alpha + \beta \in L$. Finally,
since $\alpha\beta$ and $\alpha + \beta$ each lie in $L$, we see that
$|L[\alpha]:L| \le 2$, and thus $\alpha$ and $\beta$ lie in the real
repeated radical extension $L[\alpha]$ of $Q$.\qed
\medskip

Finally, we want to show that the hypothesis in Theorem~B that $Q$
is real cannot be removed. 
\medskip

\iitem{(9.6) EXAMPLE.}~~{\sl There exist fields $Q \sbs K \sbs L \sbs \C$,
where $L$ is a repeated radical extension of $Q$ and $|L:Q|$ is odd, but
where $K$ is not a repeated radical extension of $Q$.}
\medskip

\iitem{Proof.}~~Let $L = \Q[\epsilon]$, where $\epsilon$ is a primitive
complex $19\,$th root of unity, so that $|L:\Q| = 18$. Let $Q$ be the
unique quadratic extension of $\Q$ in $L$ and let $K$ be the unique field
of degree $3$ over $Q$ in $L$. Now $L = Q[\epsilon]$ is a radical
extension of $Q$ of degree $9$ and we claim that the cubic extension
$Q \sbs K$ is not a repeated radical extension. It suffices, of course, to
show that $K$ is not a radical extension of $Q$.

The only roots of unity in $L$ are the $38\,$th roots of unity, and thus
the only roots of unity in $K$ are $\pm1$. We know that $K$
is Galois over $Q$ since $\gal LQ$ is abelian, and since $|K:Q| = 3$,
we see by Lemma~2.3(a) that $K$ cannot be radical over $Q$.\qed
\goodbreak
\bigskip\bigskip

\centerline{REFERENCES}
\parindent = 0pt
\bigskip\frenchspacing

1. I. M. Isaacs, Solutions of polynomials by real radicals, Amer. Math.
Monthly, {\bf 92} (1985) 571--575.
\smallskip

2. I. M. Isaacs, {\it Algebra: A graduate course}, Brooks/Cole,
Pacific Grove, 1994.

\bye